\input amstex
\documentstyle{amsppt}
\magnification=\magstep1

\pageheight{9.0truein}
\pagewidth{6.5truein}

\NoBlackBoxes

\input xy
\xyoption{curve}\xyoption{arrow}\xyoption{matrix}

\long\def\ignore#1{#1}

\def\sN{{\Bbb N}}

\def\pinf{{{\Cal P}^\infty}}
\def\pdim{\operatorname{p\,dim}}
\def\la{{\Lambda}}
\def\lamod{\Lambda\operatorname{-mod}}
\def\deltamod{\Delta\operatorname{-mod}}
\def\ximod{\Xi\operatorname{-mod}}
\def\laMod{\Lambda\operatorname{-Mod}}
\def\Hom{\operatorname{Hom}}
\def\lfindim{\operatorname{l\,fin\,dim}}
\def\lFindim{\operatorname{l\,Fin\,dim}}
\def\ker{\operatorname{ker}}
\def\length{\operatorname{length}}

\topmatter

\title Viewing finite dimensional representations through
infinite dimensional ones\endtitle
\rightheadtext{Finite versus infinite dimensional representations}

\author Dieter Happel and Birge Huisgen-Zimmermann\endauthor

\address Fakult\"at f\"ur Mathematik, Technische Universit\"at
Chemnitz-Zwickau, Postfach 964, D-09009 Chemnitz, Germany\endaddress

\email happel\@mathematik.tu-chemnitz.de\endemail

\address Department of Mathematics, University of California, Santa Barbara, CA
93106, USA\endaddress

\email birge\@math.ucsb.edu\endemail

\dedicatory Dedicated to the memory of Maurice Auslander\enddedicatory

\thanks The research of the second author was partially supported by a grant
from the National Science Foundation. Part of this work was done while the first author was visiting UCSB. He would like to thank his coauthor for her kind hospitality.\endthanks

\abstract We develop criteria for deciding the contravariant finiteness status of a subcategory $A \subseteq \lamod$, where $\Lambda$ is a finite dimensional algebra. In particular, given a finite dimensional $\Lambda$-module $X$, we introduce a certain class of modules -- we call them $A$-phantoms of $X$ -- which indicate whether or not $X$ has a right $A$-approximation: We prove that $X$ fails to have such an approximation if and only if $X$ has infinite-dimensional $A$-phantoms. Moreover, we demonstrate that large phantoms encode a great deal of additional information about $X$ and $A$ and that they are highly accessible, due to the fact that the class of all $A$-phantoms of $X$ is closed under subfactors and direct limits. \endabstract

\endtopmatter

\document

\head 1. Introduction and preliminaries\endhead

Given a finite dimensional algebra $\la$ and a resolving
contravariantly finite subcategory $\frak A$ of the category $\lamod$ of
finitely generated left $\la$-modules, the minimal right
$\frak A$-approximations of the simple left $\la$-modules provide significant
structural information about arbitrary objects of $\frak A$:  Indeed, if these
approximations are labeled $A_1, \dots, A_n$, then a finitely generated
$\la$-module $M$ belongs to $\frak A$ if and only if $M$ is a direct summand
of a module that has a filtration with consecutive factors in $\{A_1, \dots,
A_n\}$;  this was proved by Auslander and Reiten [2].  (For the basic
definitions, see under `Preliminaries' below.)  Among other consequences, this
result has an obvious homological pay-off: Namely, if all of the
$A_i$ have finite projective dimensions, then the maximum of these dimensions
coincides with the supremum of the projective dimensions attained on objects
of $\frak A$.

Our main interest will be in the category $\frak A=\pinf(\lamod)$ of finitely
generated
$\la$-modules of finite projective dimension, even though many of our results
will address the situation of an arbitrary subcategory of $\lamod$ which is
closed under finite direct sums.  Loosely speaking, the category 
$\pinf(\lamod)$ will be contravariantly finite in
$\lamod$ if it is either very large or very small; for example, by [3], every
representation-finite subcategory of $\lamod$ is contravariantly finite.  There
are several
 classes of algebras for which the problem of whether $\pinf(\lamod)$
is contravariantly finite in $\lamod$ is settled, e.g., for left serial
algebras, it is settled in the positive [4]. However, in general, it is
difficult to decide whether
$\pinf(\lamod)$ for a given algebra
$\la$ has this property.  In particular, there is only a single instance for
which failure of contravariant finiteness of $\pinf(\lamod)$ has actually been
established;  this is an example of a monomial relation algebra, due to Igusa,
Smal\o, and Todorov [9], which is so closely related to the Kronecker algebra
that a proof for failure of contravariant finiteness of $\pinf(\lamod)$ can be
gleaned from the representation theory of this latter algebra (of course, for
the Kronecker algebra itself, as for any algebra $\la$ of finite global
dimension, we know that $\pinf(\lamod) = \lamod$ is contravariantly finite). 
The present research was triggered by a question of M. Auslander as to
whether, for the monomial relation algebra with differing big and little
finitistic dimensions which was exhibited by the second author in [8],
$\pinf(\lamod)$ is contravariantly finite.  This specific problem turned out
to be rather intractable without a systematic theory providing direction.  

Our goal here is twofold.  On one hand, we establish `negative criteria' for
the contravariant finiteness of subcategories $\frak A$ of $\lamod$.  These
criteria are particularly manageable in the case of monomial relation
algebras, and one of them yields a negative answer to the specific
question mentioned above.  On the other hand -- maybe more importantly --  our
investigation reveals that, for non-contravariantly finite subcategories
$\frak A$ of
$\lamod$, certain modules of infinite dimension over the base field take over
the role which is played by the minimal right $\frak A$-approximations in
the contravariantly finite case.  We call these modules `$\frak A$-phantoms'
and define them as follows:  Given a subclass
$\frak C$ of $\frak A$ and a finitely generated left $\la$-module $X$, we
label as
$\frak C$-approximation of $X$ inside
$\frak A$ any homomorphism $f : A \rightarrow X$ with $A \in \frak A$ having
the property that each homomorphism $g : C \rightarrow X$ with $C \in \frak C$
factors through $f$.  Moreover, a $\la$-module $H$, not necessarily in $\frak
A$, is called an $\frak A$-phantom of $X$ in case there exists a non-empty
finite subclass $\frak C$ of
$\frak A$ such that each $\frak C$-approximation of $X$ inside $\frak A$ has
$H$ as a subfactor; direct limits of such modules $H$ are again called $\frak
A$-phantoms of $X$. In case $X$ has a right $\frak A$-approximation, the
$\frak A$-phantoms of $X$  are clearly just the subfactors of the 
minimal right $\frak A$-approximation of
$X$.  For an intuitive idea of the information stored in  such phantoms,
suppose for the moment that
$X=S$ is simple; comparable to the minimal right
$\frak A$-approximation of $S$ (in case of existence), the $\frak A$-phantoms
of
$S$ represent, in as highly compressed a form as possible, the relations
characterizing those objects of $\frak A$ which carry $S$ in their tops. 
In essence, our negative criteria give instructions for the uncovering
of phantoms which are too big (namely, non-finitely generated) to be
compatible with contravariant finiteness of
$\frak A$.  Our main result (Theorem 9) ensures existence in general: 
Namely, we prove that, for any subcategory $\frak A$ of $\lamod$ which is
closed under finite direct sums, a finitely generated $\la$-module $X$ fails
to have a right
$\frak A$-approximation if and only if $X$ has $\frak A$-phantoms of infinite
dimension over the base field.
	
The core of the paper is Section 3. The concepts introduced at the beginning
of that section appear cogent in the light of the preliminary criterion for
failure of contravariant finiteness presented in Section 2 and its
applications to special cases.  We present several
examples to illustrate how phantoms mark the dividing line between
contravariant finiteness and failure thereof, and to show concretely what
type of information they store.  These examples also indicate how sensitive
the property of contravariant finiteness of 
$\pinf(\lamod)$ in $\lamod$ is, with respect to
modifications of the relations of the underlying algebra. \medskip

\noindent {\it Acknowledgment:\/} The authors would like to thank Sverre
Smal\o{} for his careful reading of a preliminary manuscript and several
helpful suggestions.
\medskip

\noindent {\it Prerequisites:} Throughout, $\la$ will denote a basic finite
dimensional algebra over a field $K$, with a fixed set of primitive
idempotents, $e_1,\dots,e_n$, and $J$ will be the Jacobson radical of $\la$.
The simple left $\la$-modules $\la e_i/Je_i$ will be abbreviated by
$S_i$. Moreover, $\lamod$ will stand for the category of all finitely
generated left $\la$-modules, and $\laMod$ for the category of {\it all} left
$\la$-modules.

Let $\frak A$ be a full subcategory of $\lamod$. Following
Auslander/Smal\o\ [3] and Auslander/Reiten [2], we say that a module
$M\in
\lamod$ has a (right) {\it $\frak A$-approximation} in case there exists a
homomorphism $\varphi : A\rightarrow M$ with $A\in \frak A$ such that each
homomorphism $B\rightarrow M$ with $B\in \frak A$ factors through $\varphi$.
By [2], existence of any $\frak A$-approximation of $M$ entails existence of
a {\it minimal} $\frak A$-approximation of $M$, i.e., one of minimal
$K$-dimension, which is unique up to isomorphism. If each object in $\lamod$
has an $\frak A$-approximation, then $\frak A$ is said to be {\it
contravariantly finite in $\lamod$} (see [2,3]). Our favorite choice of a
subcategory $\frak A\subseteq \lamod$ will be the subcategory of all finitely
generated left $\la$-modules of finite projective dimension; we label it
$\pinf(\lamod)$.  Finally, we call a module category $\frak A \subseteq
\lamod$ {\it resolving} in case $\frak A$ is closed under extensions, as well
as kernels of epimorphisms, and contains all indecomposable projective left
$\la$-modules.  Clearly, $\pinf(\lamod)$ is an instance of a resolving
subcategory of $\lamod$.  

Given a module $M\in \laMod$, we call an element $m\in M$ a {\it top element}
of $M$ if $m\in M\setminus JM$ and $e_im=m$ for some $i$; in this case, we
also say that $m$ is a top element {\it of type $e_i$}. In all of our
examples, $\la$ will be a split finite dimensional algebra over $K$, that is,
$\la$ will be of the form $K\Gamma/I$ where $\Gamma$ is a quiver and $I$ an
admissible ideal in the path algebra $K\Gamma$. We will briefly and
informally review the second author's conventions for the graphical
communication of information about countably generated $\la$-modules. For
additional detail, see [7,8]. (We point out that our labeled graphs are related
to the module diagrams studied by Alperin [1] and Fuller [5].) 

Let $\Gamma$ be the quiver 

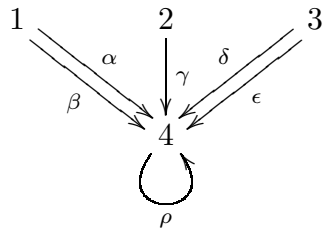
\ignore{
$$\xymatrixcolsep{3pc}
\xymatrix{
1 \ar[dr]<0.5ex>^\alpha \ar[dr]<-0.5ex>_\beta &2 \ar[d]^\gamma &3
\ar[dl]<-0.5ex>_\delta \ar[dl]<0.5ex>^\epsilon\\
 &4 \ar@'{@+{[0,0]+(-5,-5)}@+{[0,0]+(0,-10)}@+{[0,0]+(5,-5)}}_\rho
}$$
} 

\noindent and $\la= K\Gamma/\langle \rho^2\rangle$. To say that a left
$\la$-module $M$ has the {\it layered and labeled graph} shown below, 

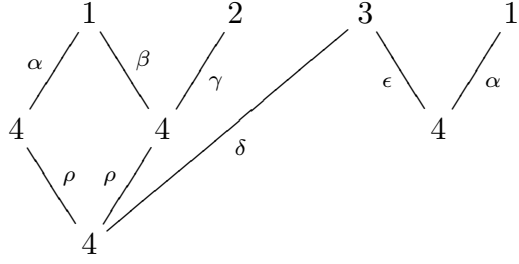
\ignore{
$$\xymatrixcolsep{1pc}
\xymatrix{
 &1 \ar@{-}[dl]_\alpha \ar@{-}[dr]^\beta &&2 \ar@{-}[dl]^\gamma &&3
\ar@{-}[ddllll]^\delta \ar@{-}[dr]_\epsilon &&1 \ar@{-}[dl]^\alpha\\
4 \ar@{-}[dr]^\rho &&4 \ar@{-}[dl]_\rho &&&&4\\
 &4
}$$
}

\noindent with
respect to a sequence $m_1,m_2,m_3,m_4$ of top elements of $M$ which are
$K$-linearly independent modulo $JM$, is to convey the following information:

$\bullet$ $M/JM\cong S_1^2\oplus S_2\oplus S_3$, the top elements $m_i$  of
$M$ have type $e_i$ for $i=1,2,3$, and $m_4$ has
type
$e_1$;

$\bullet$ $JM/J^2M\cong S_4^3$,  the three copies of $S_4$ modulo $J^2M$
being generated by $\alpha m_1$, $\beta m_1$, $\epsilon m_3$, and such that
$\gamma m_2$ is congruent to a nonzero scalar multiple of $\beta m_1$ modulo
$J^2M$, and $\alpha m_4$ congruent to a nonzero scalar multiple of $\epsilon
m_3$;

$\bullet$ $J^2M/J^3M= J^2M\cong S_4$ is generated by $\rho\alpha m_1$, and
is also generated by any of the elements $\rho\beta m_1$, $\rho\gamma m_2$,
or $\delta m_3$.

Finally, we will consider the following two finitistic dimensions of $\la$:
the {\it left little finitistic dimension}, $\lfindim \la$, which is the
supremum of the finite projective dimensions attained on $\lamod$, and the
{\it left big finitistic dimension}, $\lFindim \la$, which stands for the
analogous supremum attained on all of $\laMod$.

\head 2. A few motivating examples\endhead

In this section, we specialize to the situation where the algebra $\la$ is 
split, i.e., we assume throughout that $\la = K \Gamma / I$ is a path algebra
modulo relations. We start by exhibiting an elementary sufficient condition for
failure of contravariant finiteness of $\pinf (\lamod)$ in $\lamod$. In our
first example we will apply it to the algebra constructed by
Igusa/Smal\o/To\-dor\-ov [9], which provided the first known instance of
such failure. While this criterion is quite easy to handle, its scope is
rather limited, and it will later be supplemented by a criterion of far
wider applicability.

For convenience of exposition, we will often view left $\la$-modules $M$ as
representations of the quiver $\Gamma$ satisfying the commutativity
relations dictated by $I$; given a path $p : e_1\rightarrow e_2$ in
$K\Gamma$, we will in that case, write $f_p : e_1M\rightarrow e_2M$ for the
$K$-linear map corresponding to $p$.

\proclaim{Elementary Criterion 1} Let $\la = K \Gamma/I$.
Suppose that $e_1$ and $e_2$ are vertices of the quiver $\Gamma$ (not
necessarily distinct) and $p, q \in K \Gamma \setminus I$ paths from $e_1$ to
$e_2$ with $\la p \cap \la q = 0$ (we view $p$ and $q$ as
elements of $\la$ whenever indicated by the context).  Moreover,
suppose that

{\rm (i)} the cyclic module
$\la (p,q)$ generated by the element 
$(p,q) \in \la^2$ has finite projective dimension, 

\noindent and that one of the
following two conditions is satisfied:  either,

{\rm (ii)} whenever $M \in \pinf (\lamod)$, then $f_p(e_1M\setminus JM) \cap
f_q(e_1JM) =\varnothing$;

\noindent or,

{\rm (ii')} whenever $M \in \pinf (\lamod)$, then $\ker(f_p)\subseteq
\ker(f_q)$, and $\ker(f_p) \subseteq e_1JM$. 

\noindent Then the simple module $S_1 =
\la e_1/J e_1$ does not have a right $\pinf 
(\lamod)$-ap\-prox\-i\-ma\-tion.\endproclaim

Before we justify the criterion, we point out that
 the second part of Criterion 1(ii') can often be verified without effort;
we label it as follows:

(iii) Whenever $M$ in $\pinf
(\lamod)$, 
then $\ker(f_p)\subseteq e_1JM$.

\noindent Indeed, suppose that $p$ is an arrow such that $\la p$ splits off
in $Je_1$ (this is obviously true when $\la$ is a monomial relation algebra).
Then Condition (iii) holds if and only if $\pdim \la p =\infty$. For, if
$\pdim \la p =\infty$ and $M\in\lamod$ contains a top element of type $e_1$
which is annihilated by $p$, then $\la p$ is isomorphic to a direct summand
of $\Omega^1(M)$, which entails that $\pdim M=\infty$. If, on the other hand,
$\pdim \la p <\infty$, then the module $M=\la/\la p$ violates Condition
(iii).

A readily recognizable situation in which the blanket hypothesis of
the criterion, as well as conditions (i) and (ii) are satisfied is as
follows: $p$ is an arrow $e_1 \rightarrow e_2$,
$q \in K\Gamma \setminus I$ a path from $e_1$ to $e_2$ of
positive length which is different from $p$ such that
$Jp = qJ = 0$, and $\pdim \la q < \infty$, while $\pdim \la e_2/J e_2 =
\infty$.

\demo{Proof of Criterion 1} We start by assuming the blanket
hypothesis of the criterion and condition (i) to construct an infinite
family of objects
$(N_n)_{n\in
\Bbb N}$ in $\pinf(\lamod)$.  Namely, for $n \in \Bbb N$, we let $b_1 =
\dots  = b_n = e_1$, define a left $\la$-module 
$$N_n : =\left( \bigoplus^n_{i=1} \la b_i \right) \bigg/ \left(
\sum^{n-1}_{i=1}
\la (p b_i - q b_{i+1}) \right),$$
and write $\overline{b_i}$ for
the residue class of $b_i$ in $N_n$. To compute the first syzygy $\Omega^1
(N_n)$ of $N_n$, consider the projective cover $\pi:
\bigoplus^n_{i=1} \la b_i \rightarrow N_n$ with $\pi (b_i) =
\overline{b_i}$, and set $C_i = \la (p
b_i - q b_{i+1})$ for
$1 \le i \le n-1$.  Note that $C_i \simeq \la (p, q)$, whence
$\pdim C_i < \infty$ by condition (i). We will conclude that $\pdim N_n
< \infty$ by showing that  $\ker \pi = \bigoplus^{n-1}_{i=1} C_i$.
Suppose that 
$$0 = \sum^{n-1}_{i=1} \lambda_i (p b_i - q b_{i+1}) =
\lambda_1 p b_1 + (\lambda_2 p - \lambda_1 q) b_2 + \dots + (\lambda_{n-1}
p - \lambda_{n-2} q) b_{n-1} - \lambda_{n-1} q b_n$$
for certain coefficients $\lambda _i \in \la$.
This implies that $\lambda_1 p = 0$, $\lambda_{n-1} q = 0$ and $\lambda_{i}
p - \lambda_{i-1} q = 0$ for $2 \le i \le n-1$, and in view of the hypothesis
that
$\la p \cap \la q = 0$, the latter equations entail $\lambda_i p =
\lambda_{i-1} q = 0$.  Thus we obtain $\lambda_i (p b_i - q b_{i+1}) = 0$ for
$1 \le i \le n-1$ as required.

Case I.  Suppose that, in addition, condition (ii) holds, but that
 there nonetheless exists a right $\pinf(\lamod)$-approximation $\varphi : A
\rightarrow S_1$ for $S_1$.  Pick $n >\length (A)$, define
$f: N_n \rightarrow S_1$ via
$f(\overline{b_1}) = e_1 + J e_1$, $f (\overline{b_i}) = 0$ for $i \ge 2$,
and let $g \in  \Hom_\Lambda (N_n, A)$ be such that $f = \varphi g$, i.e., such
that the following diagram commutes:
  
\ignore{
$$\xy\xymatrix{  &N_n\ar[dl]_g\ar[d]^f\\
A\ar[r]^\varphi&S_1}\endxy$$
}

Next pick $m \ge 1$ minimal with the property that $g (\overline{b_1}), \dots,
g (\overline{b_m})$ are $K$-linearly dependent modulo $J A$; such an integer
$m$ exists because $\length (N_n/J N_n) = n >\length (A)$. Say $\sum^m_
{i=1} k_i g (\overline{b_i}) \in J A$ with $k_i \in K$, not all zero. Clearly,
$k_m \ne 0$.  Moreover, $k_1 = 0$, since $0 = f (\sum^m_{i=2}
k_i
\overline{b_i}) = \varphi g (\sum^m_{i=2} k_i \overline{b_i}) = - k_1 \varphi
g (\overline{b_1}) = - k_1 (e_1 + J e_1)$. In particular, this shows that
$m\ge 2$. Now set $x = \sum^m_{i=2} k_i \overline b_{i-1}$. We will check that
$g (x) =
\sum^{m-1}_{i=1} k_{i+1} g (\overline{b_i})$ again belongs to $J A$,
a contradiction to the minimal choice of $m$. Indeed, $p g (x) = g (
\sum^m_{i=2} k_i p \overline{b}_{i-1}) = g (\sum^m_{i=2} k_i q \overline {b_i}) = q g
(\sum^m_{i=2} k_i \overline{b_i}) \in q J A$, which by condition (ii) entails
that $g (x)$ is not a top element of $A$. This shows $g(x)\in JA$ and thus
completes the argument for Case I.

Case II. Now suppose that conditions (i) and (ii') hold. Again assume that $
S_1$ has a right $\pinf (\lamod)$-ap\-prox\-i\-ma\-tion $\varphi: A \rightarrow
S_1$, and for $n >\length (A)$, define $f: N_n \rightarrow S_1$ as in case I.
In turn choose $g \in  \Hom_\Lambda (N_n, A)$ such that $f = \varphi g$.  But
this time, pick $m \in \sN$ minimal with the property that $\sum^ m_{i=1} k_i g
(\overline{b_i}) \in J A$ and $q (\sum^m_{i=1} k_i g (\overline {b_i})) = 0$,
for some scalars $k_i \in K$ which are not all zero. Such an $m$ exists because
$\ker (g)$ intersects the $K$-subspace of $N_n$ generated by $\overline
b_1,\dots, \overline b_n$ non-trivially. As before we obtain $k_1 = 0$ and $m
\ge 2$, and again, we set
$x =
\sum^m_{i=2} k_i \overline {b_{i-1}}$ and compute $p g (x) = q (\sum^m_{i=2}
k_i g (\overline{b_i}))=0$.  Now condition (ii') guarantees that $g (x) =
\sum^{m-1}_{i=1} k_ {i+1} g (\overline{b_i})$ is not a top element of $A$ and
that $q g (x) = q (\sum^{m-1}_{i=1} k_{i+1} g (\overline{b_i})) = 0$. This
is, once more, incompatible with the minimal choice of $m$. \qed\enddemo

\example{Example 2} [9]
Let $\Lambda = K \Gamma/I$ be the monomial relation algebra with quiver
  
\ignore{
$$\xy\xymatrix{\Gamma:&1\ar@/^.3pc/[r]^\beta\ar@/^1.2pc/[r]^\alpha
&2\ar@/^.2pc/[l]^\gamma}\endxy$$
}
  
\noindent and ideal $I = \langle \alpha\gamma, \beta \gamma, \gamma
\beta\rangle$ of relations.  Then the indecomposable projective left
$\Lambda$-modules have the following graphs:
  
\ignore{
$$\xy\xymatrixcolsep{.7pc}\xymatrix{
&1\ar@{-}[ld]_\alpha \ar@{-}[rd]^\beta&&&&2\ar@{-}[d]_\gamma\\
2\ar@{-}[d]_\gamma&&2&&&1\\
1
}\endxy$$
}

Set $p = \beta$ and $q = \alpha$, and verify the simplified versions of
conditions (i) and (ii), as spelled out before the proof of Criterion 1:
Clearly, $\alpha J =J \beta =  0$ and $\pdim \la \alpha = 0 < \infty$,
while $\pdim\Lambda e_2/Je_2 = \infty$.  Thus the criterion guarantees that
$S_1$ does not have a right 
$\pinf (\lamod)$-ap\-prox\-i\-ma\-tion. \qed\endexample

Note that a typical class of modules defeating attempts to find a
right $\pinf (\lamod)$-ap\-prox\-i\-ma\-tion of $S_1$ in Example 2 is
the following class $(M_n)_{n \in \sN}$ of strings of composition length $2
n$, uniquely determined by their graphs:
  
\ignore{
$$\xy\xymatrixcolsep{.7pc}\xymatrix{
\save +<-1pc,-1.7pc>\drop{M_n:}&1\ar@{-}[dr]_(0.42)\beta
&&1\ar@{-}[dl]_(0.58)\alpha
\ar@{-}[dr]_(0.42)\beta&&\save+<0pc,-1.7pc>\drop{\cdots}
&&1\ar@{-}[dl]_(0.58)\alpha\ar@{-}[dr]_(0.42)\beta\\ &&2&&2&&2&&2
}\endxy$$
}
  
\noindent (Observe that this family of modules $M_n$ represents just a minor
simplification of the family of test modules $N_n$ exhibited in the proof of
Criterion 1 for the choices
$p=\beta$ and
$q=\alpha$; namely $M_n \cong N_n / \la \alpha \overline {b_1}$.)  A stumbling
block for contravariant finiteness can be more succinctly communicated via the
infinite dimensional direct limit of the strings $M_n$:
  
\ignore{
$$\xy\xymatrixcolsep{.7pc}\xymatrix{
\save+<-2pc,-1.7pc>\drop{\varinjlim M_n:}&1\ar@{-}[dr]_(0.42)\beta
&&1\ar@{-}[dl]_(0.58)\alpha
\ar@{-}[dr]_(0.42)\beta&&1\ar@{-}[dr]_(0.42)\beta\ar@{-}[dl]_(0.58)\alpha
&&\save+<0pc,-1.7pc>\drop{\cdots}\\ &&2&&2&&2
}\endxy$$
}    
  
\noindent At the same time, this limit is a `minimal' module 
$M \in \pinf (\laMod)$ with the property 
that all homomorphisms $M_n \rightarrow S_1$, $n \in \sN$, factor through
$M$.

By the preceding example one could be led to believe that -- under the
assumption that $\la \beta$ is a direct summand of $Je_1$ having infinite
projective dimension -- the existence of a family of modules
$(M_n)_{n \in
\Bbb N}$ as above should imply non-existence of a $\pinf$-approximation of the
simple module associated with the vertex `$1$'.  This is not the case,
however.  In fact, the existence of modules $(M_n)_{n \in \sN}$ of the
indicated shape  inside
$\pinf (\lamod)$ is only potentially troublesome. Continuing to denote the
algebra of Example 2 by $\la$, we next give
an example of an algebra $\Delta$ with $\lamod \subseteq\deltamod$
such that $\Delta e_i = \la e_i$ for $i = 1,2$ and $(\lamod) \cap
\pinf(\deltamod) =
\pinf(\lamod)$; in particular, all of the
$\Lambda$-modules $M_n$ belong also to $\pinf(\deltamod)$. On the other
hand, we will see that $\pinf (\deltamod)$ {\it is} contravariantly finite in
$\deltamod$ in this example.

\example{Example 3} Let $\Delta = K  \Gamma'/I'$, where
$\Gamma'$ is   
  
\ignore{
$$\xy\xymatrix{3\ar[r]^\delta&1\ar@/^.3pc/[r]^\beta\ar@/^1.2pc/[r]^\alpha
&2\ar@/^.2pc/[l]^\gamma}\endxy$$
}
  
\noindent and the ideal $I'$ is generated by monomial relations such that the
indecomposable projective left $\Delta$-modules have graphs
  
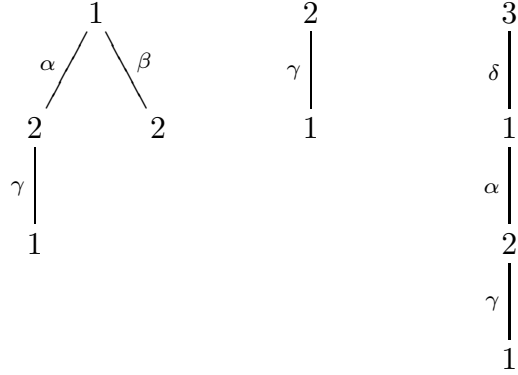
\ignore{
$$\xy\xymatrixcolsep{.7pc}\xymatrix{
&1\ar@{-}[ld]_\alpha \ar@{-}[rd]^\beta&&&&2\ar@{-}[d]_\gamma
&&&&3\ar@{-}[d]_\delta\\ 2\ar@{-}[d]_\gamma&&2&&&1&&&&1\ar@{-}[d]_\alpha\\
1&&&&&&&&&2\ar@{-}[d]_\gamma\\
&&&&&&&&&1
}\endxy$$
} 
  
\noindent This time the `zipper effect' of the previous example (where the 
requirement that the homomorphisms $f \in \Hom_\Lambda (M_n, S)$
be factorizable through a module 
$M \in \pinf (\lamod)$
forces the $K$-dimension of $M$ to grow with increasing $n$) can be stopped.
This is due to the new projective module $\Lambda e_3$.

We will prove the contravariant finiteness of 
$\pinf (\deltamod)$ by
explicitly describing minimal 
$\pinf (\deltamod)$-ap\-prox\-i\-ma\-tions of
the simple left $\Delta$-modules. We claim that the following canonical
epimorphism is a right $\pinf (\deltamod)$-ap\-prox\-i\-ma\-tion of $S_1$:
  
\ignore{
$$\xy\xymatrixcolsep{.7pc}\xymatrix{
&&&1\ar@{-}[ddr]_\beta &&3\ar@{-}[d]^\delta\\
\varphi_1:&A_1&=&&&1 \ar@{-}[ld]^\alpha&\ar[rr]&&&S_1\\
&&&&2
}\endxy$$
}
  
\noindent Note that $A_1 = \bigl(\Delta e_1 \oplus \Delta e_3 \bigr)/
\bigl( \Delta \alpha e_1 + \Delta (\beta e_1 -
\alpha \delta e_3) \bigr)$ is the injective envelope of the
simple  module $S_2= \Delta e_2/Je_2$.  We again use
$J$ to denote the Jacobson radical of
$\Delta$. and if $x_1$ and $x_3$ stand for the residue classes of
$e_1$ and $e_3$ in $A_1$, respectively, we let $\varphi_1 (x_1) = e_1 + J
e_1$ and $\varphi_1 (x_3) = 0$.

Note first that $\Omega^1 (A_1) = (\Delta e_2)^2$, whence
$\pdim A_1 = 1$.  Now let $f : M \rightarrow S_1$ be an epimorphism with $M
\in \pinf (\deltamod)$, $m_1 \in M$ a top element of type $e_1$ such
that $f(m_1) = e_1 + Je_1$, and let $m_2,\dots,m_r\in \ker(f)$ be such that
$m_1 + JM, \dots,  m_r + JM$ form a $K$-basis for $e_1(M/JM)$. Then
each
$\beta m_i$ is a nonzero element in the socle of $M$, since otherwise
$\Omega^1(M)$ would contain a direct summand isomorphic to $S_2$, which is
impossible in view of the fact that
$\pdim S_2 = \infty$, while $\pdim M <\infty$. More strongly, this argument
shows that $\beta m_1,\dots, \beta m_r$ are $K$-linearly independent, whence
in particular $\beta \overline m_1$ remains nonzero in
$\overline M = M / \left(\sum_{i=2}^2 \Delta \beta m_i \right)$.  It is 
enough to factor the map $\overline f : \overline M \rightarrow S_1$
induced by $f$ through $\varphi_1$.  In doing this, it is clearly harmless to
assume that
$\alpha \overline m_1= k\beta \overline m_1$ for some scalar $k\in K$ which
may be zero; if this is not a priori the case, we factor out the submodule
$\Delta (\alpha -\beta) \overline m_1$ in addition.  Furthermore, we may
assume that the elements $m_2, \dots , m_r \in \ker(f)$ are chosen in such a
way that there exists an integer $s$ between $2$ and $r$ with the property
that $\beta \overline m_1 = \alpha \overline m_i$ for $ 2 \le i \le s$ and
$\Delta \beta \overline m_1 \cap \left( \sum_{i=s+1}^r \Delta \alpha
\overline m_i \right) =0$.

Set $B=
\sum_{i=1}^r \Delta \overline m_i$, and let $\iota : B\rightarrow \overline
M$ be the canonical embedding.  In view of the prededing adjustments, we can
define a map
$\sigma \in \Hom_\Delta(B,A_1)$ by setting $\sigma(\overline m_1) = x_1 +
k \delta x_3$, $\sigma(\overline m_i)=\delta x_3$ for $ 2 \le i \le s$, and
$\sigma(\overline m_i)=0$ for $s+1 \le i \le r$.
Since $A_1$ is injective, $\sigma$ can be extended to a homomorphism $\tau
\in \Hom_\Delta(\overline M,A_1)$ which makes the lower triangle in the
diagram below commute.

\ignore{
$$\xy\xymatrixcolsep{3pc}\xymatrix{
A_1 \ar[r]^{\varphi_1} &S_1\\
B \ar[u]^\sigma \ar[r]^\iota &{\overline M} \ar[u]_{\overline f}
\ar@{-->}[ul]_\tau }\endxy$$
}

\noindent Our construction entails that the upper triangle then commutes as
well, which shows that $\varphi_1$ is indeed a (minimal)
$\pinf(\deltamod)$-approximation of $S_1$.

It is less involved to see that the canonical maps

\ignore{
$$\xy\xymatrixcolsep{.2pc}\xymatrixrowsep{.3pc}\xymatrix{
&&2\ar@{-}[dd]_\gamma&&&&&&&&&3\ar@{-}[dd]_\delta\\
\varphi_2:&&&\ar[rrr]&&&&S_2&\text{\quad and \quad}&\varphi_3:&&&\ar[rrr]
&&&&S_3\\ &&1&&&&&&&&&1
}\endxy$$
}

\noindent are right 
$\pinf (\deltamod)$-ap\-prox\-i\-ma\-tions of $S_2$ and $S_3$, respectively.
By [2], this shows that $\pinf (\deltamod)$ is contravariantly finite in
$\deltamod$.\qed\endexample

Example 3 also shows that the hypotheses of Criterion 1
cannot be simplified to the combination of conditions (i) and (iii), where
(iii) is as in the remark following the statement of the criterion.
Indeed, if in Example 3, we take $p = \beta$ and $q = \alpha$,
then both (i) and (iii) are satisfied, but $S_1$ does admit a right 
$\pinf$-ap\-prox\-i\-ma\-tion.

Finally, we modify the algebra $\Delta$ of Example 3 very slightly to  an
algebra $\Xi$, with the effect that $\pinf (\ximod)$ again
fails to be contravariantly finite. Here condition (ii) fails for any
choice of $p$ and $q$, but (i) and (ii') are satisfied.  This sequence of
modifications illustrates a phenomenon which will become more obvious in the
sequel: Namely, that contravariant finiteness of
$\pinf(\lamod)$ in $\lamod$ --  as well as failure of this condition -- is
highly unstable.

\example{Example 4} The quiver of $\Xi$ is that of the algebra
$\Delta$ in Example 3, but we delete one of the relations, with the effect
that the
$K$-dimension of $\Xi$ exceeds that of $\Delta$ by 1, and the
indecomposable projective left
$\Xi$-modules take on the form
  
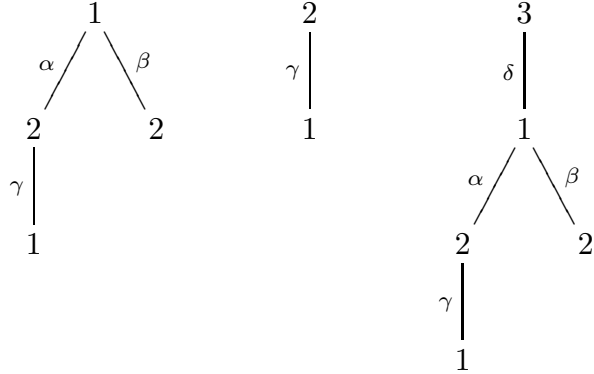
\ignore{
$$\xy\xymatrixcolsep{.7pc}
\xymatrix{
&1\ar@{-}[ld]_\alpha \ar@{-}[rd]^\beta&&&&2\ar@{-}[d]_\gamma
&&&&3\ar@{-}[d]_\delta\\
2\ar@{-}[d]_\gamma&&2&&&1&&&&1\ar@{-}[dl]_\alpha\ar@{-}[dr]^\beta\\
1&&&&&&&&2\ar@{-}[d]_\gamma&&2\\ &&&&&&&&1
}\endxy
$$
}
  
\noindent We again apply Criterion 1 with the choice $p = \beta$ and $q =
\alpha$. Clearly, $\Lambda \alpha \cap \Lambda \beta = 0$. Moreover,
$\pdim\Xi (\beta, \alpha) = \pdim \left( \ignore{ \xymatrixrowsep{.5pc}
 \xy<1pc,1pc>
\xymatrix{2\ar@{-}[d]^\gamma\\ 1 }\endxy} \right) = 0$,
whence condition (i) of the criterion is satisfied.

Next, we check that condition (ii') of our criterion is satisfied. Let
$M$ be in $\pinf (\ximod)$. As in Example 3, it is readily 
checked that, given any top element $x \in M$ of type $e_1$, we have $\beta
x \ne 0$. So, in proving that for an {\it arbitrary} element $x \in M$, the
vanishing of $\beta x$ implies the vanishing of $\alpha x$, we may assume
that $x = e_1x \in J M$ with $\beta x = 0$; since $\beta J = K \beta
\delta$ and $\alpha J = K \alpha \delta$, we may moreover assume that $x \in 
\delta M$. Suppose that $\alpha x \ne 0$.  In view of the equality $\alpha
\delta J = 0$, this implies that $x \in \delta M \backslash \delta J M$, i.e.,
$x = \delta y$ for some top element $y \in M$ of type  $e_3$.
Let $\pi: P =
\Xi e_3 \oplus Q \rightarrow M$ be a projective cover with $\pi (e_3)
= y$. Then $\beta \delta e_3$ is a nonzero element of $\ker \pi = \Omega^1
(M)$, and  since $\Xi \beta \delta e_3 \simeq S_2$ has infinite projective
dimension and is thus not a direct summand of $\Omega^1 (M)$, we see that
$\beta \delta e_3 \in J\Omega^1 (M) \cap \beta JP = \beta \Omega^1 (M)$;
the last equality follows from [5, Lemma 1]. Thus $\beta \delta e_3 =
\beta z$, where $z = e_1 z \in e_1\Omega^1 (M)e_3 \subseteq e_1J Pe_3$.  
Since, clearly, the desired implication `$\beta u = 0 \implies \alpha u =
0$' does hold for arbitrary elements $u$ of a projective left
$\Xi$-module, we deduce that $\alpha\delta e_3 =\alpha z$ and
conclude that
$\alpha x =
\alpha \pi (z) = 0$ as required.

Finally, we note that condition (ii) fails in this example. Indeed, if it
would hold, it would be true for $p=\beta$ and $q=\alpha$.  However,the left
$\Xi$-module
$$M = \bigl( \Xi e_1 \oplus \Xi e_3 \bigr)/\bigl( \Xi \alpha e_1 +
\Xi (\beta e_1 - \alpha \delta e_3) \bigr)$$
with graph
  
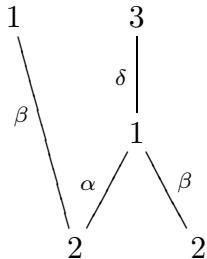
\ignore{
$$\xy\xymatrixcolsep{.7pc}\xymatrix{
1\ar@{-}[ddr]_(.4)\beta &&3\ar@{-}[d]_\delta\\
&&1\ar@{-}[dl]_\alpha\ar@{-}[dr]^\beta\\
&2&&2
}\endxy
$$
}
  
\noindent has syzygy $\Omega^1 (M) = \ignore{ \xymatrixrowsep{.5pc}
 \xy<1pc,1pc>
\xymatrix{2\ar@{-}[d]^\gamma\\ 1 }\endxy \oplus
\xy<1pc,1pc>
\xymatrix{2\ar@{-}[d]^\gamma\\ 1 }\endxy }$,
and thus $\pdim M = 1$, but if $x_i$ is the residue class of $e_i$ in
$M$ for $i = 1, 3$, then $x_1$ is a top element of type $e_1$ with $\beta
x_1 = \alpha \delta x_3 \in \alpha J M$. \qed\endexample

As in Example 2, we can again -- in the preceding example -- pin down classes of
modules of finite projective dimension which are responsible for failure of
contravariant finiteness of
$\pinf (\ximod)$. For instance, there is no homomorphism $\varphi : A
\rightarrow S_1$ with $A \in \pinf (\ximod)$ such that {\it all} the canonical
epimorphisms from the modules
  
\ignore{
$$\xy\xymatrixcolsep{.7pc}\xymatrix{
&&1\ar@{-}[ddr]_\beta
&&3\ar@{-}[d]_\delta&&3\ar@{-}[d]_\delta&&&&3\ar@{-}[d]_\delta\\
E_n:&&&&1\ar@{-}[dl]_(0.55)\alpha
\ar@{-}[dr]_(0.42)\beta&&1\ar@{-}[dl]_(0.55)\alpha
\ar@{-}[dr]_(0.42)\beta\ar@{{}{}}[rrrr]|{\txt{$\cdots$}}
&&&&1\ar@{-}[dl]_(0.55)\alpha\ar@{-}[dr]_(0.42)\beta\\ &&&2&&2&&2&&2&&2
}\endxy
$$
}
  
\noindent onto $S_1$ can be factored through $\varphi$. Observe, however, that
they can all be factored through the canonical surjection from
  
\ignore{
$$\xy\xymatrixcolsep{.7pc}\xymatrix{
&&1\ar@{-}[ddr]_\beta
&&3\ar@{-}[d]_\delta&&3\ar@{-}[d]_\delta&&3\ar@{-}[d]_\delta\\ E=\varinjlim
E_n:&&&&1\ar@{-}[dl]_(0.55)\alpha
\ar@{-}[dr]_(0.42)\beta&&1\ar@{-}[dl]_(0.55)\alpha
\ar@{-}[dr]_(0.42)\beta&&1\ar@{-}[dl]_(0.55)\alpha\ar@{-}[dr]_(0.42)\beta
&&\cdots\\ &&&2&&2&&2&&2
}\endxy
$$
}
  
\noindent onto $S_1$.

\head 3. Relative approximations and phantoms\endhead

If ${\frak A} \subset \lamod$ is a resolving contravariantly finite
subcategory of $\lamod$, then, according to [2], the minimal right
approximations of the simple left $\Lambda$-modules hold a substantial
amount of information on arbitrary objects of ${\frak A}$. The gap in the
available information when ${\frak A}$ is not contravariantly finite is to be
filled by direct limits of `partial approximations' as indicated informally
in Section 2.  (We follow the
convention that `direct limits' are colimits extending over directed index
sets.)

\definition{Definitions 5} Let ${\frak C} \subset {\frak A}$ be full 
subcategories of $\lamod$ such that ${\frak A}$ is closed under finite direct
sums, and let $\widehat{\frak A}$ be the closure of ${\frak A}$ under direct limits
in
$\laMod$. Moreover, let $X$ be a finitely generated left 
$\Lambda$-module.

(1) A {\it (right) ${\frak C}$-ap\-prox\-i\-ma\-tion of $X$ inside ${\frak
A}$} (resp\. {\it inside $\widehat{\frak A}$}) is a homomorphism $f:A
\rightarrow X$ with
$A$ in ${\frak A}$ (resp\. $A$ in $\widehat{\frak A}$) such that 
$$\Hom(-,A)|_{\frak C} @>\Hom(-,f)>> \Hom(-,X)|_{\frak C} @>>> 0$$
is an exact sequence of functors, i.e., such that each map $g
\in \Hom_\Lambda (C, X)$ with $C$ in ${\frak C}$ factors through $f$.

If ${\frak C} = {\frak A}$, a ${\frak C}$-ap\-prox\-i\-ma\-tion of $X$
inside ${\frak A}$ will simply be called a (right) {\it ${\frak
A}$-ap\-prox\-i\-ma\-tion of} $X$, in accordance with the existing terminology.

(2) A {\it ${\frak C}$-phantom of $X$ relative to ${\frak A}$ of the
first kind} is an object $B$ in $\lamod$ (not necessarily in ${\frak A}$)
with the following property: There exists a finite non-empty set ${\frak C}
(B) \subset {\frak C}$ such that, for each ${\frak C}
(B)$-ap\-prox\-i\-ma\-tion 
$f:A
\rightarrow X$  inside ${\frak A}$, the module $B$ is a subfactor of $A$.
Any direct limit of ${\frak C}$-phantoms of $X$ of the first kind will be
called a {\it ${\frak C}$-phantom of $X$ of the second kind}. 

We will refer to both kinds of phantoms as {\it ${\frak C}$-phantoms} of
$X$ relative to
${\frak A}$ and, more briefly, to {\it ${\frak A}$-phantoms} if ${\frak C} =
{\frak A}$.

(3) A ${\frak C}$-phantom $B$ of $X$ relative to ${\frak A}$ is
called {\it effective} if there exists a homomorphism $f: B \rightarrow X$
which is a ${\frak C}$-ap\-prox\-i\-ma\-tion of $X$ inside $\widehat{\frak A}$ (in
particular, $B \in \widehat{\frak A}$ in that case). \enddefinition
 \medskip

In case $X$ fails to have an $\frak A$-approximation, there may be a plethora
of $\frak A$-phantoms. This facilitates the construction of infinite
dimensional phantoms, which in turn lie at the heart of our criteria for
failure of contravariant finiteness of $\frak A$.

The {\it effective} phantoms, on the other hand, are in a sense the best
possible substitutes for minimal approximations in the sense of
Auslander/Smal\o\ [3] and Auslander/Reiten [2]. For example, if $X = \Lambda
e/J e$ is simple, the effective $\frak A$-phantoms of $X$ compress information
about the relations of those modules in
${\frak A}$ which have a top element of type $e$ into the tightest possible
format. Compared with classical approximations, we simply renounce the
requirement that this picture  should fit into a finitely generated module.

Our primary interest here will be in the situation where 
${\frak A} = \pinf
(\lamod)$ and ${\frak C}$ is a countable family of modules of finite
projective dimension.

\example{Basic Observations 6} Let ${\frak C}, {\frak A}, \widehat{\frak A}$
and $X$ be as in Definition 5.

(1) Since $\widehat{\frak A}$ is closed under arbitrary direct sums, it is
clear that
$X$ always has ${\frak C}$-ap\-prox\-i\-ma\-tions inside $\widehat{\frak
A}$. 
Just 
add up a sufficient number of copies of each object $C$ in ${\frak C}$ to
cover all homomorphisms $C \rightarrow X$.  Similarly, $X$ has a ${\frak
C}$-ap\-prox\-i\-ma\-tion inside ${\frak A}$  whenever
${\frak C} \subseteq {\frak A}$ is a finite subset, because 
$\Hom_\Lambda (C, X)$
has finite $K$-dimension for each $C \in {\frak C}$. 

(2) If $X$ has phantoms relative to
${\frak A}$ which have unbounded finite lengths, then $X$ fails to have an
${\frak A}$-ap\-prox\-i\-ma\-tion. In particular, this is true in case $X$ has
a non-finitely generated ${\frak A}$-phantom.

(3) Suppose that, in addition to being closed under finite direct sums,
${\frak A}$
is closed under direct summands. Then the existence of an ${\frak
A}$-ap\-prox\-i\-ma\-tion of $X$ implies the existence of a {\it unique}
minimal such  approximation by [2], say $A(X)$. In that case, $A(X)$ is the
only effective ${\frak A}$-phantom
of $X$, and all other ${\frak A}$-phantoms of $X$ are subfactors of $A(X)$.

(4)  If ${\frak C} \subseteq {\frak C}' \subseteq {\frak A}$, each ${\frak
C}$-phantom relative to ${\frak A}$ is also a ${\frak C}'$-phantom relative
to ${\frak A}$.  So, in particular, each $\frak C$-phantom relative to $\frak
A$ is an $\frak A$-phantom. 
\qed\endexample

\example{7. The Examples of Section 2 Revisited} Let $X = S_1 = \Lambda e_1 /
J e_1$, and 
${\frak A} = \pinf (\lamod)$.

In Example 2, the module
  
\ignore{
$$\xy\xymatrixcolsep{.7pc}\xymatrix{
\save+<-3pc,-1.7pc>\drop{M=\varinjlim M_n:}&1\ar@{-}[dr]_(0.42)\beta
&&1\ar@{-}[dl]_(0.58)\alpha
\ar@{-}[dr]_(0.42)\beta&&1\ar@{-}[dr]_(0.42)\beta\ar@{-}[dl]_(0.58)\alpha
&&\save+<0pc,-1.7pc>\drop{\cdots}\\ &&2&&2&&2
}\endxy
$$
}
  
\noindent is an effective ${\frak C}$-phantom of $S_1$ inside 
$\pinf (\laMod)$, where ${\frak C} = \{M_n \mid n \in \sN\}$, but
  
\ignore{
$$\xy\xymatrixcolsep{.7pc}\xymatrix{
&&1\ar@{-}[dr]_(0.42)\beta\ar@{-}[dl]_(0.58)\alpha&&1\ar@{-}[dl]_(0.58)\alpha
\ar@{-}[dr]_(0.42)\beta&&1\ar@{-}[dr]_(0.42)\beta\ar@{-}[dl]_(0.58)\alpha
&&\save+<0pc,-1.7pc>\drop{\cdots}\\
\save+<-1.5pc,1pc>\drop{N:}&2\ar@{-}[d]_\gamma&&2&&2&&2\\
&1
}\endxy
$$
}
  
\noindent is neither the source of a ${\frak C}$-ap\-prox\-i\-ma\-tion 
of $S_1$ inside
$\pinf (\laMod)$, nor a ${\frak C}$-phantom of $S_1$ relative to
$\pinf(\lamod)$.

As for Example 3:  Start by observing that the above graphs uniquely define
left modules over the modified algebra
$\Delta$, again denoted $M_n, M$ and $N$, and the class of $\Delta$-modules
${\frak C} = \{ M_n \mid n \in \Bbb N \}$ in turn belongs to
$\pinf(\deltamod)$.  Moreover, the homomorphism
$f: M\rightarrow S_1 = \Delta e_1/J e_1$ which sends the top element
represented by the left-most `$1$' in the graph of $M$ to $e_1 + Je_1$ and
sends the top elements displayed farther to the right to zero is still
a ${\frak C}$-ap\-prox\-i\-ma\-tion of $S_1$ inside
$\pinf (\deltamod)$.  However, in the present setup, both $M$ and $N$ fail
to be  ${\frak C}$-phantoms of $S_1$ relative to 
$\pinf(\deltamod)$; indeed, as we saw earlier, $S_1$ has a 
$\pinf (\deltamod)$-ap\-prox\-i\-ma\-tion in that example.

Finally, let us focus on Example 4. Viewing $M_n$ and $M$ as left
$\Xi$-modules, and keeping in mind that the $M_n$ belong to $\pinf
(\ximod)$, we find that $M$ is an effective ${\frak C}$-phantom of $S_1$
relative to $\pinf (\ximod)$,  as in Example 2.  Moreover, if ${\frak
E} = \{E_n \mid n \in
\Bbb N \}$ with $E_n$ as defined after Example 4, then $M$ is also an 
${\frak E}$-phantom of $S_1$ relative to $\pinf
(\ximod)$, 
but not an effective one because the canonical epimorphism
  
\ignore{
$$\xy\xymatrixcolsep{.7pc}\xymatrix{
&&1\ar@{-}[ddr]_\beta &&3\ar@{-}[d]_\delta\\
E_2&=&&&1 \ar@{-}[ld]^\alpha&\ar[rr]&&&S_1\\
&&&2
}\endxy$$
}
  
\noindent does not factor through $M$. The better ${\frak E}$-phantom here is 
$E = \varinjlim E_n$,
which is actually an effective $({\frak C} \cup {\frak E})$-phantom of
$S_1$ relative to $\pinf (\ximod)$. \qed\endexample

Next we prepare for a general existence result.  In a nutshell: Whenever
${\frak A} \subseteq
\lamod$ fails to be contravariantly finite, there exist ${\frak A}$-phantoms of
infinite $K$-dimension.

\proclaim{Proposition 8} Suppose that ${\frak A} \subseteq \lamod$
is closed under finite direct sums and that $X \in \lamod$ does not
have a (right) ${\frak A}$-ap\-prox\-i\-ma\-tion. Then there exists a countable
subclass ${\frak C}$ of ${\frak A}$ such that $X$ fails to have a 
${\frak C}$-ap\-prox\-i\-ma\-tion inside ${\frak A}$.\endproclaim

\demo{Proof} By repeatedly applying the first of the observations under 6, 
we show that, for each $d \ge 1$, there exists a finite subset ${\frak C}_d$
of ${\frak A}$ such that $X$ does not have a ${\frak C}_d$-ap\-prox\-i\-ma\-tion of
$K$-dimension $\le d$ inside ${\frak A}$.

Assuming the contrary for some $d \ge 1$, 
we pick a module $Y_1$ in ${\frak A}$
with $\Hom_\Lambda (Y_1, X) \ne 0$  -- such a module $Y_1$ exists by
hypothesis  --  and let $f_1: A_1 \rightarrow X$ be a 
$\{Y_1\}$-ap\-prox\-i\-ma\-tion of $X$ inside ${\frak A}$ such that dim$_k A_1 \le d$.
In particular, $f_1$ is nonzero. Since $X$ fails to have an 
${\frak A}$-ap\-prox\-i\-ma\-tion, there exists an object $Y_2$ in ${\frak A}$
such that some homomorphism in $\Hom_\Lambda (Y_2, X)$ fails to factor through
$f_1$. Let $f_2: A_2 \rightarrow X$ be an $\{A_1, Y_2\}$-ap\-prox\-i\-ma\-tion of
$X$ inside ${\frak A}$; by assumption, $A_2$ can be chosen to have 
$K$-dimension
at most $d$. Inductively, our assumption thus yields a family $(A_n)_{n \ge 1}$
of objects of ${\frak A}$ with dim$_k A_n \le d$ for all $n$, together with
finitely generated left $\Lambda$-modules $Y_n$ and homomorphisms $f_n:
A_n \rightarrow X$ such that $f_n$ is an $\{A_{n-1},
Y_n \}$-ap\-prox\-i\-ma\-tion of
$X$ inside ${\frak A}$, but fails to be a $\{Y_{n+1}\}$-ap\-prox\-i\-ma\-tion.
Accordingly, we can pick $g_n \in \Hom_\Lambda (A_n, A_{n+1})$ such
that $f_n = f_{n+1}g_n$ and none of the $g_n$ is an isomorphism. But since
$f_1 = f_2 g_1 = f_3 g_2 g_1 = \dots = f_{n+1} g_n \cdots g_1$ is nonzero,
we deduce $g_n \cdots g_1 \ne 0$ for all $n$, which contradicts the
Harada-Sai Lemma [6] and proves our assumption to be absurd.

Letting ${\frak C}_d$ for $d\ge 1$ be as in our initial claim, the
countable subset ${\frak C} = \bigcup_{d \ge 1} {\frak C}_d$
of ${\frak A}$ is clearly as desired. \qed\enddemo

We apply this proposition to obtain the announced existence result.

\proclaim{Theorem 9} Suppose that ${\frak A} \subseteq \lamod$ is
closed under finite direct sums, and let $X \in \lamod$. Then the 
following conditions are equivalent:

{\rm (1)} $X$ fails to have an ${\frak A}$-ap\-prox\-i\-ma\-tion.

{\rm (2)} There exists a countable subclass ${\frak C} \subseteq {\frak A}$
such that $X$ has an effective ${\frak C}$-phantom of countably infinite
$K$-dimension relative to ${\frak A}$.

{\rm (3)} $X$ has an ${\frak A}$-phantom of infinite
$K$-dimension. \endproclaim

\demo{Proof}
`(1)$\implies$(2)'. Assume that (1) holds.  Then Proposition 8 yields a
countable subclass
${\frak D} = \{D_1, D_2, D_3, \dots \}$ of ${\frak A}$ such that $X$ does
not have a ${\frak D}$-ap\-prox\-i\-ma\-tion inside ${\frak A}$.  However, by
the first of the observations under 6, there exists a 
$\{D_1 \}$-ap\-prox\-i\-ma\-tion of $X$ inside ${\frak A}$, say $f_1: A_1
\rightarrow X$.

Next we pick an $\{A_1\}$-ap\-prox\-i\-ma\-tion $f_2: A_2 \rightarrow X$ of $X$
inside ${\frak A}$, together with a map $g_{1,2} \in \Hom_\Lambda
(A_1, A_2)$ satisfying $f_1 = f_2 \circ g_{1,2}$, 
such that dim$_K \bigl(g_{1,2} (A_1) \bigr)$ is as small
as possible. Consequently, the following is true: Whenever $f'_2: A'_2
\rightarrow X$ is an $\{A_2\}$-ap\-prox\-i\-ma\-tion of $X$ inside ${\frak A}$
and $g' \in \Hom_\Lambda (A_2, A'_2)$ is such that $f_2 = f'_2 \circ
g'$, we have $g' (g_{1,2} (A_1)) \simeq g_{1,2} (A_1)$.

We now choose any $\{D_2, A_2\}$-ap\-prox\-i\-ma\-tion $f_3: A_3 \rightarrow X$
of
$X$ inside ${\frak A}$, and subsequently an $\{A_3\}$-ap\-prox\-i\-ma\-tion 
$f_4: A_4 \rightarrow X$ inside ${\frak A}$, together with a map $g_{3,4} \in
\Hom_\Lambda (A_3, A_4)$ such that $f_3 = f_4 \circ g_{3,4}$ and
dim$_K \bigl( g_{3,4} (A_3) \bigr)$ is minimal.

Continuing along this line, we obtain a sequence of objects $(A_n)_{n \ge
1}$ and maps $f_n: A_n \rightarrow X$ such that, for $n \ge 2, \, f_{2n-1}$
is a $\{D_n, A_{2n-2}\}$-ap\-prox\-i\-ma\-tion of $X$ inside ${\frak A}$ and
$f_{2n}$ is an $\{A_{2n-1}\}$-ap\-prox\-i\-ma\-tion which is coupled with a
map $g_{2n-1, 2n} \in \Hom_\Lambda (A_{2n-1}, A_{2n})$ such that
$f_{2n-1} = f_{2n} \circ g_{2n-1, 2n}$ and dim$_K \bigl(g_{2n-1,2n}
(A_{2n-1}) \bigr)$ is minimal.

Set ${\frak C} = \{A_1, A_2, A_3, \dots \}$ and supplement the above maps 
$g_{n,n+1}$
for odd $n$ by homomorphisms $g_{n,n+1} \in \Hom_\Lambda (A_n, A_{n+1})$
with $f_n = f_{n+1} \circ g_{n,n+1}$ for $n$ even. If, for $n < m$, we moreover
define $g_{n,m} = g_{m-1,m} \circ \dots \circ g_{n,n+1}: A_n \rightarrow
A_m$, then $(A_n, g_{n,m})_{n, m \in \sN, n < m}$ is an inductive system with
$f_n = f_m \circ g_{n,m}$. Set 
$$A = \varinjlim A_n, \qquad \qquad f=\varinjlim 
f_n \in \Hom_\Lambda (A,X),$$ 
and let $h_n: A_n \rightarrow A$ be the canonical maps. Clearly, $A$ belongs to
$\widehat{\frak A}$ (see Definition 5).  Moreover, each homomorphism in
$\Hom_\Lambda (C, X)$ with $C \in {\frak C}$ factors through $f$ and, a
fortiori,  so does each homomorphism
in $\Hom_\Lambda (D_n, X)$. In other words, $f : A \rightarrow X$ is a
${\frak C} \cup {\frak D}$-approximation of $X$ inside $\widehat{\frak A}$.

Next we want to identify $A$ as a ${\frak C}$-phantom of $X$ relative 
to ${\frak A}$.  Our construction entails that, for
$m > 2n$, we have $g_{2n, m} \circ g_{2n-1, 2n} (A_{2n-1}) \simeq g_{2n-1, 2n}
(A_{2n-1})$, and consequently we have $h_{2n} (U_{2n}) \simeq U_{2n}$ if we
define 
$U_{2n}= g_{2n-1,2n} (A_{2n-1})$.
Since $A$ is the directed union of the submodules $h_{2n} (U_{2n}), \,
n \in\sN$, it suffices to show that each of the modules $U_{2n}$ is a 
${\frak C}$-phantom of $X$ relative to ${\frak A}$ of the first kind.   
For that purpose, consider the finite subset ${\frak C} (U_{2n}) =
\{A_{2n}\}$ of ${\frak C}$, and let $f': A' \rightarrow X$ be a ${\frak C}
(U_{2n})$-ap\-prox\-i\-ma\-tion of $X$ inside ${\frak A}$. If $g' \in
\Hom_\Lambda (A_{2n}, A')$ is such that $f_{2n} = f' \circ g'$, our
construction yields $U_{2n}
\simeq g' (U_{2n}) \subseteq A'$, which shows that $U_{2n}$ is indeed a
${\frak C}$-phantom of $X$ relative to ${\frak A}$ of the first kind.
Consequently, $A$ is a 
${\frak C}$-phantom of $X$ relative to ${\frak A}$ of the second kind which,
by the preceding paragraph, is even effective.

Finally, we note that dim$_K A \le \aleph_0$ by construction.  To prove the
reverse inequality, we assume, to the contrary, that dim$_k A <
\infty$.
But this means that $f$ is a ${\frak D}$-ap\-prox\-i\-ma\-tion of
$X$ inside ${\frak A}$, which contradicts our choice of ${\frak D}$ and
completes the proof of `(1)$\implies$(2)'.

The implications `(2)$\implies$(3)' and `(3)$\implies$(1)' are immediate
consequences of the basic observations 6(4) and 6(2), respectively. 
\qed\enddemo

The following is an upgraded version of the elementary Criterion 1 for
non-existence of a $\pinf(\lamod)$-ap\-prox\-i\-ma\-tion of a given simple
module $S$. The idea underlying the proof is the same, even though we impose
no restrictions on the subcategory $\frak A \subseteq \lamod$ this time. In
particular, this criterion again points to a countable subclass 
${\frak C}$ of $\frak A$ which obstructs the approximability of
$S$ by a {\it finitely generated} module of finite projective dimension.  In
view of the proof of Theorem 9, it can hence be used towards the explicit
construction of $\frak A$-phantoms of $S$. 
While this criterion will be instrumental in resolving the problem of
contravariant finiteness in our key example (Section 4), for complex
non-monomial algebras, it may still be nontrivial to verify or refute
Condition (2) below.  We therefore add an illlustration of how the underlying
idea can still be used towards deciding questions of contravariant finiteness,
even when the criterion is not readily applicable verbatim.

\proclaim{Criterion 10} Suppose that $\la$ is a split finite dimensional
algebra and
$\frak A$ a full subcategory of
$\lamod$.  Moreover,
let $e_1, \dots, e_m$ be pairwise  orthogonal primitive idempotents of
$\Lambda$, and $p_1,
\dots, p_m, q_1, \dots, q_m \in J$ with $p_i = p_i e_i$ and $q_i = q_i e_i$
such that the following conditions are satisfied:

{\rm (1)} For each $n \in \sN$, 
there is a module $M_n \in \frak A$, together with a sequence $x_{n1}, \dots,
x_{n,mn}$ of $mn$ top elements of $M_n$ which are $K$-linearly independent
modulo
$J M_n$ such that $0 \ne p_{r(i)} x_{ni} = q_{r(i+1)} x_{n,i+1}$ for $1 \le i
< mn$, where $r (i) \in \{1, \dots, m\}$ is congruent to $i$ modulo $m$.

{\rm (2)} For any object $C$ in $\frak A$, the following are true:

\quad{\rm (i)} if $x \in C$ is a top element of type $e_1$ 
then $p_1 x \ne 0$;

\quad{\rm (ii)}  if $y, z \in C$ with $0 \ne p_{r(i)} y =
q_{r(i+1)} z$, 
then $p_{r(i+1)} z \ne 0$.

Then $S_1 = \Lambda e_1/J e_1$ does not have an 
$\{ M_n \mid  n \in \Bbb N \}$-ap\-prox\-i\-ma\-tion inside $\frak A$.  In
particular, $\frak A$ is not contravariantly finite in $\lamod$ in
that case.\endproclaim

\demo{Proof}  Assume to the contrary that $f: A \rightarrow S_1$
is an $\{ M_n \mid  n \in \Bbb N \}$-ap\-prox\-i\-ma\-tion of $S_1$, and
choose
$n \in \Bbb N$ such that dim$_K A < mn$.  Fixing $n$,  we will briefly
write $x_i$ for $x_{ni}$, $1 \le i \le mn$.  Consider the homomorphism
$g: M_n
\rightarrow S_1$, defined by $g (x_1) = e_1 + J e_1$ and $g (x_i) = 0$ for $2
\le i \le mn$; this definition is meaningful because $x_1$ is a top element of
$M_n$ of type $e_1$. Choose $h: M_n 
\rightarrow A$ such that $g = f h$. Moreover, note that, due to the linear
dependence of the elements $h (x_1), \dots, h (x_{mn})$ of $A$, there
exists a natural number $t$, together with scalars $k_t, \dots, k_n$ such
that $k_t \ne 0$ whereas
$0 = p_{r (t)}
\sum^{mn} _{i=t} k_i h (x_i) = p_{r (t)} z$, where $z = \sum^{mn}_
{i=t} h (e_{r(t)} k_i x_i)$.

Choose $t$ minimal with the property that the next-to-last equation is
satisfied for some nonzero scalar $k_t$ and some scalars $k_{t+1}, \dots ,
k_{mn}$. We claim that
$t > 1$. Indeed, if
$t = 1$, then $f(z) = g (\sum e_1 k_i x_i) = k_1 g (x_1) \ne 0$ and hence $z$
is a top element of
$A$ of type $e_1$. By (2)(i) this implies that $p_1 z \ne 0$, a contradiction.
Now set $y = \sum^{mn}_{i=t} h (e_{r(t-1)} k_i x_{i-1})$. By the minimal
choice of $t$, we obtain $p_{r(t-1)} y \ne 0$. Using condition (1) and the
choice of the $r(i)$, we further compute that $p_{r(t-1)} y = q_{r(t)} z$, and
 -- invoking condition (2)(ii) -- we conclude $p_{r(t)} z \ne 0$. But this is
incompatible with our choice of
$t$. \qed\enddemo

As is backed up by the proof of Criterion 10, the hypotheses of this
criterion entail the existence of an infinite dimensional phantom of
$S_1$ relative to $\frak A$, 
a graph of which contains a subgraph
  
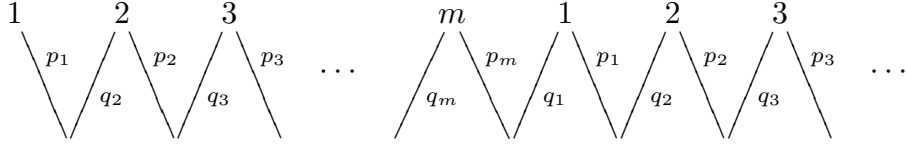
\ignore{
$$\xy\xymatrixcolsep{.7pc}\xymatrix{
1\ar@{-}[dr]^(0.40){p_1}&&2\ar@{-}[dl]^(0.60){q_2}
\ar@{-}[dr]^(0.40){p_2}&&3\ar@{-}[dr]^(0.40){p_3}\ar@{-}[dl]^(0.60){q_3}
&&*!/u1.5pc/{\cdots}&&m\ar@{-}[dr]^(0.40){p_m}\ar@{-}[dl]^(.60){q_m}
&&1\ar@{-}[dl]^(0.60){q_1}\ar@{-}[dr]^(0.40){p_1}&&2\ar@{-}[dl]^(0.60){q_2}
\ar@{-}[dr]^(0.40){p_2}&&3\ar@{-}[dr]^(0.40){p_3}\ar@{-}[dl]^(0.6){q_3}
&&*!/u1.5pc/{\cdots}\\ &&&&&&&&&&&&&&&&&
}\endxy$$
}
  
\noindent As mentioned earlier, the criterion works well for
$\frak A = \pinf (\lamod)$ when $\la$ is a monomial relation algebra.  One of
the reasons for this can be found in the following observation which shows how
easy it is to get a chain of
$\pinf( \lamod)$-phantoms started in the monomial situation.

\example{Remarks 11} Let $\Lambda= K\Gamma/I$ be a monomial relation algebra,
and suppose that the simple module $S_1 = \Lambda e_1 / J e_1$ has infinite
projective dimension. If $\alpha_1, \dots, \alpha_r$ are arrows $\alpha_j:
e_1 \rightarrow e_j$ ending in distinct vertices $e_1, \dots, e_r$, such
that $\pdim \la \alpha_j = \infty$ for $1 \le j \le r$, then
  
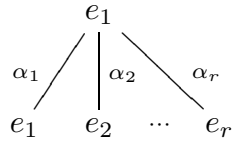
\ignore{
$$\xy\xymatrixcolsep{.7pc}\xymatrix{
&e_1\ar@{-}[dl]_(0.65){\alpha_1}\ar@{-}[d]^(0.57){\alpha_2}
\ar@{-}[drr]^(0.7){\alpha_r}\\ e_1&e_2\ar@{{}{}}[rr]|\cdots&&e_r
}\endxy$$
}
  
\noindent is a subgraph of the graph of a $\pinf (\lamod)$-phantom $C$ of
$S_1$.  More precisely, there exists a top element $c\in C$ of type $e_1$
such that $\alpha_i c \ne 0$ for $1 \le i \le r$. 

If, moreover, there exists a module 
$M \in \pinf (\lamod)$ with a graph containing a subgraph of the form
  
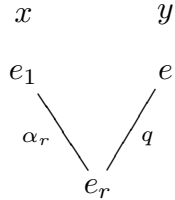
\ignore{
$$\xy\xymatrixcolsep{.7pc}\xymatrix{
e_1 \ar@{-}[dr]_{\alpha_r} \save+<0pc,1.5pc>\drop{x} &&e \ar@{-}[dl]^q
\save+<0pc,1.5pc>\drop{y}\\
 &e_r
}\endxy$$
}
  
\noindent where $x, y \in M$ are top elements of types $e_1$ and $e$
respectively and $q$ denotes a path in $K\Gamma \setminus I$, then there
exists a 
$\pinf (\lamod)$-phantom of $S_1$ whose graph contains a subgraph of the form

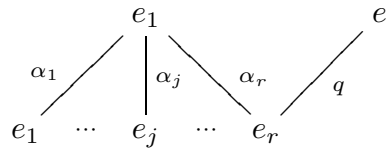
\ignore{
$$\xy\xymatrixcolsep{.7pc}\xymatrix{
 &&e_1\ar@{-}[dll]_(0.65){\alpha_1}\ar@{-}[d]^(0.57){\alpha_j}
\ar@{-}[drr]^(0.7){\alpha_r} &&&&e \ar@{-}[dll]^q\\ 
e_1\ar@{{}{}}[rr]|\cdots &&e_j\ar@{{}{}}[rr]|\cdots &&e_r
}\endxy$$
}

\noindent with respect to suitable top elements.
\endexample

\demo{Proof} The second statement clearly follows from the
first. To justify the first, we start by noting that $\bigoplus_
{1 \le j \le r} \Lambda \alpha_j$ is a direct summand of $J e_1$, due to
the fact that $\Lambda$ is a monomial relation algebra. 
Let $C \in \pinf
(\lamod)$ have a top element $c$ of type $e_1$. To see that $\alpha_j c \ne 0$
for $1 \le j \le r$, consider a projective cover
$$\pi: P = \Lambda x_0 \oplus \bigoplus_{i \in I} \Lambda x_i 
\rightarrow C$$
such that $x_0 = e_1$ and $\pi (x_0) = c$. If we had $\alpha_j c = 0$,
we could conclude that $\Lambda \alpha_j x_0$ is a direct summand of
$\ker \pi = \Omega^1 (C) \le J P$, because $\Lambda \alpha_j x_0$ is a
direct summand of $J P$, which is incompatible with our setup.
\qed\enddemo

While most of our applications demonstrate the use of phantoms towards a
proof that $\pinf (\lamod)$ fails to be contravariantly 
finite, phantoms may also be helpful in  finding 
$\pinf (\lamod)$-ap\-prox\-i\-ma\-tions.

\example{Example 12} [2, Example on p. 137] Let $\Lambda = K \Gamma/I$
be based on the quiver
  
\ignore{
$$\xy\xymatrix{
1\ar@'{@+(-5,5)@+(-10,0)@+(-5,-5)}_\alpha\ar[r]^\beta&2\ar[r]^\gamma
&3\ar[r]^\delta&4 }\endxy$$
}
  
\noindent such that the $\Lambda e_i$ have graphs:
  
\ignore{
$$\xy\xymatrixcolsep{.7pc}
\xymatrix{
&1\ar@{-}[ld]_\alpha \ar@{-}[rd]^\beta&&&&2\ar@{-}[d]_\gamma
&&&&3\ar@{-}[d]_\delta&&&&4\\
1\ar@{-}[d]_\beta&&2\ar@{-}[d]^\gamma&&&3\ar@{-}[d]_\delta&&&&4\\
2\ar@{-}[d]_\gamma&&3&&&4\\ 3
}\endxy$$
}
  
\noindent Clearly, $S_1$ is the only simple left $\Lambda$-module of infinite
projective dimension. By Remark 11, $S_1$ has a 
$\pinf (\lamod)$-phantom $A_1$,
with graph 
$\ignore{  \xymatrixrowsep{.5pc}
\xy<1pc,1pc>
\xymatrix{1\ar@{-}[d]^\alpha\\
1
}\endxy }$,
and it is readily checked that there is no object in $\pinf(\lamod)$ having a
submodule with graph 
$\ignore{ 

\xymatrixrowsep{.5pc}\xymatrixcolsep{0.5pc}
\xy<1pc,1pc>
\xymatrix{1 \ar@{-}[dr]_\alpha &&e \ar@{-}[dl]\\
 &1
}\endxy }$.
 One deduces that $A_1\rightarrow S_1$ is a (minimal)
$\pinf (\lamod)$-ap\-prox\-i\-ma\-tion
of $S_1$. \qed\endexample

For non-monomial relation algebras, one often needs to slightly vary
the idea of Criterion 10. We illustrate the construction of phantoms in
such a non-monomial situation.

\example{Example 13} Let $\Lambda = K \Gamma /I$, where $\Gamma$ is
the quiver
  
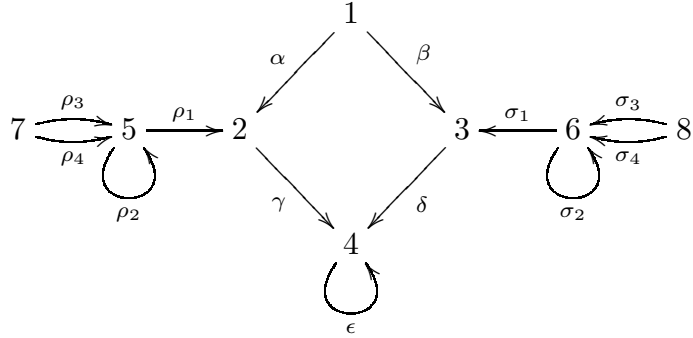
\ignore{
$$\xy\xymatrix{
&&&1\ar[dl]_\alpha\ar[dr]^\beta\\
7\ar@/^.3pc/[r]^{\rho_3}\ar@/_.3pc/[r]_{\rho_4}
&5\ar@'{@+{[0,0]+(-5,-5)}@+{[0,0]+(0,-10)}@+{[0,0]+(5,-5)}}_{\rho_2}
\ar[r]^{\rho_1}&2\ar[dr]_\gamma&&3\ar[dl]^\delta
&6\ar[l]_{\sigma_1}
\ar@'{@+{[0,0]+(-5,-5)}@+{[0,0]+(0,-10)}@+{[0,0]+(5,-5)}}_{\sigma_2}
&8\ar@/^.3pc/[l]^{\sigma_4}\ar@/_.3pc/[l]_{\sigma_3}\\
&&&4\ar@'{@+{[0,0]+(-5,-5)}@+{[0,0]+(0,-10)}@+{[0,0]+(5,-5)}}_\epsilon
}\endxy$$ 
}
  
\noindent and $I \subset K \Gamma$ is the unique ideal containing $\gamma
\alpha-
\delta\beta$ and having the property that the indecomposable projective left
$\Lambda$-modules have the graphs
  
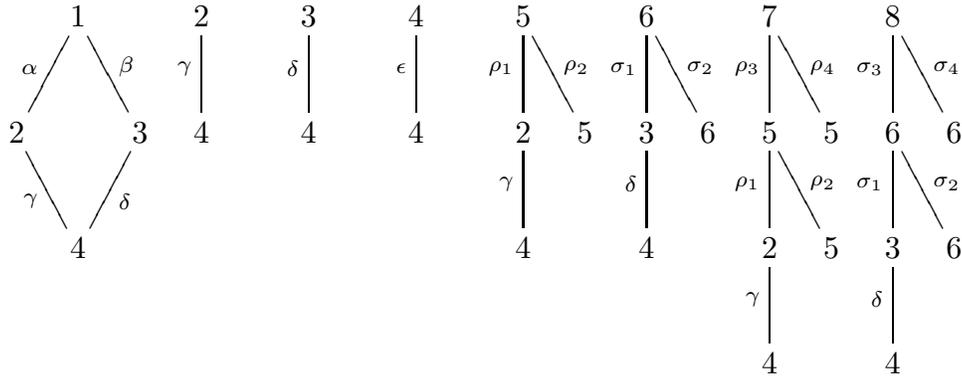
\ignore{
$$\xy\xymatrixcolsep{.7pc}
\xymatrix{
&1\ar@{-}[dl]_\alpha\ar@{-}[dr]^\beta&&2\ar@{-}[d]_(0.43)\gamma
&&3\ar@{-}[d]_(0.43)\delta&&4\ar@{-}[d]_(0.43)\epsilon
&&5\ar@{-}[d]_(0.43){\rho_1}\ar@{-}[dr]^{\rho_2}
&&6\ar@{-}[d]_(0.43){\sigma_1}\ar@{-}[dr]^{\sigma_2}
&&7\ar@{-}[d]_(0.43){\rho_3}\ar@{-}[dr]^{\rho_4}
&&8\ar@{-}[d]_(0.43){\sigma_3}\ar@{-}[dr]^{\sigma_4}\\
2\ar@{-}[dr]_\gamma&&3\ar@{-}[dl]^\delta&4&&4&&4&&2\ar@{-}[d]_(0.43)\gamma
&5&3\ar@{-}[d]_(0.43)\delta&6&5\ar@{-}[d]_(0.43){\rho_1}\ar@{-}[dr]^{\rho_2}
&5&6\ar@{-}[d]_(0.43){\sigma_1}\ar@{-}[dr]^{\sigma_2}&6\\
&4&&&&&&&&4&&4&&2\ar@{-}[d]_(0.43)\gamma&5&3\ar@{-}[d]_(0.43)\delta&6\\
&&&&&&&&&&&&&4&&4 }\endxy$$
}

We will see that $S_1 = \la e_1 / Je_1$ does not have a right $\pinf
(\lamod)$-approximation by constructing $\pinf(\lamod)$-phantoms of infinite
$K$-dimension.  We start by observing that for each module $M$ in
$\pinf (\lamod)$ with top element $m$ of type $e_1$, 
either $\alpha m \ne 0$ or $\beta m \ne 0$. Thus the module

\ignore{
$$\xy\xymatrixcolsep{.7pc}
\xymatrixrowsep{.4pc}
\xymatrix{
&1\ar@{-}[dd]&&1\ar@{-}[dd]\\
C_1=&\ar@{{}{}}[rr]|\oplus&&\\
&2&&3
}\endxy$$
}

\noindent of finite projective dimension is an (effective)
$\operatorname{add}(C_1)$-phantom of
$S_1$ inside $\pinf (\lamod)$; a fortiori, 
$C_1$ is a $\pinf
(\lamod)$-phantom of $S_1$. Since there exist modules in $\pinf(\lamod)$
having subgraphs 
  
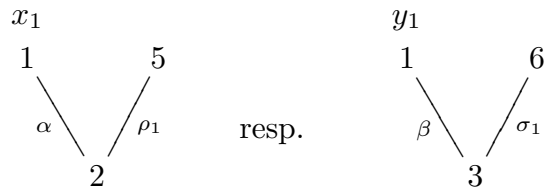
\ignore{
$$\xymatrixcolsep{.7pc}\xymatrixrowsep{.14pc}
\xy<0pc,3pc>\xymatrix{
x_1\\
1\ar@{-}[ddddr]_\alpha&&5\ar@{-}[ddddl]^{\rho_1}\\
\\
\\
\\
&2
}\endxy
\text{\qquad resp. \qquad}
\xy<0pc,3pc>\xymatrix{
y_1\\
1\ar@{-}[ddddr]_\beta&&6\ar@{-}[ddddl]^{\sigma_1}\\
\\
\\
\\
&3
}\endxy$$
}
  
\noindent where $x_1$ and $y_1$ stand for top elements, namely 
  
\ignore{
$$\xymatrixcolsep{.7pc}
\xy<0pc,2pc>\xymatrix{
1\ar@{-}[dr]_\alpha&&5\ar@{-}[dl]^{\rho_1}\ar@{-}[dr]^(.65){\rho_2}\\
&2&&5
}\endxy
\text{\qquad resp. \qquad}
\xy<0pc,2pc>\xymatrix{
1\ar@{-}[dr]_\beta&&6\ar@{-}[dl]^{\sigma_1}\ar@{-}[dr]^(.65){\sigma_2}\\
&3&&6
}\endxy$$
}
  
\noindent and since each $M \in \pinf (\lamod)$ 
with a top element $m$ of type
$e_5$ (resp\. $e_6$) satisfies $\rho_2 m \ne 0$ (resp\. $\sigma_2 m\ne
0$), the module
  
\ignore{
$$\xymatrixcolsep{.7pc}
C_2=\quad\xy<0pc,2pc>\xymatrix{
1\ar@{-}[dr]_\alpha&&5\ar@{-}[dl]^{\rho_1}\ar@{-}[dr]^(.65){\rho_2}\\
&2&&5
}\endxy
\qquad \oplus \qquad
\xy<0pc,2pc>\xymatrix{
1\ar@{-}[dr]_\beta&&6\ar@{-}[dl]^{\sigma_1}\ar@{-}[dr]^(.65){\sigma_2}\\
&3&&6
}\endxy$$
}
  
\noindent is an  $\operatorname{add}(C_2)$-phantom of $S_1$ inside $\pinf
(\lamod$).

In the next step, we observe that $\pinf (\lamod$) 
contains objects
with subgraphs
  
\ignore{
$$\xymatrixcolsep{.7pc}
\xy<0pc,2pc>\xymatrix{
1\ar@{-}[dr]_\alpha&&5\ar@{-}[dl]^{\rho_1}\ar@{-}[dr]^(.65){\rho_2}
&&7\ar@{-}[dl]^{\rho_3}\\ &2&&5
}\endxy
\text{\qquad resp. \qquad}
\xy<0pc,2pc>\xymatrix{
1\ar@{-}[dr]_\beta&&6\ar@{-}[dl]^{\sigma_1}\ar@{-}[dr]^(.65){\sigma_2}
&&8\ar@{-}[dl]^{\sigma_3}\\ &3&&6
}\endxy$$
}
  
\noindent where again $x_1$ and $x_2$, resp\. $y_1$ and $y_2$, denote top
elements, and since each module
$M \in \pinf (\lamod)$ with top element $m$ of type $e_7$ (resp\.
$e_8$) satisfies $\rho_4 m \ne 0$ (resp\. $\sigma_4 m \ne 0$) the
module
  
\ignore{
$$\xymatrixcolsep{.7pc}
C_3=\quad\xy<0pc,2pc>\xymatrix{
1\ar@{-}[dr]_\alpha&&5\ar@{-}[dl]^{\rho_1}\ar@{-}[dr]^(.65){\rho_2}
&&7\ar@{-}[dl]^{\rho_3}\ar@{-}[dr]^(.65){\rho_4}\\ &2&&5&&5
}\endxy
\qquad \oplus \qquad
\xy<0pc,2pc>\xymatrix{
1\ar@{-}[dr]_\beta&&6\ar@{-}[dl]^{\sigma_1}\ar@{-}[dr]^(.65){\sigma_2}
&&8\ar@{-}[dl]^{\sigma_3}\ar@{-}[dr]^(.65){\sigma_4}\\ &3&&6&&6
}\endxy$$
}
  
\noindent is an effective  $\operatorname{add}(C_3)$-phantom of $S_1$ 
inside $\pinf (\lamod).$
A fortiori, $C_3$ is a $\pinf (\lamod)$-phantom inside 
$\pinf(\lamod).$

Proceeding in this fashion, we obtain modules 
$C_n \in \pinf (\lamod)$ of length $4n$, namely
  
\ignore{
$$\xymatrixcolsep{.43pc}
C_n=\quad\xy<0pc,2pc>\xymatrix{
1\ar@{-}[dr]&&5\ar@{-}[dl]\ar@{-}[dr]&&7\ar@{-}[dl]\ar@{-}[dr]
&\save+<0.7pc,-1.7pc>\drop{\cdots}&&7\ar@{-}[dl]\ar@{-}[dr]\\ &2&&5&&5&5&&5
}\endxy
\quad \oplus \quad
\xy<0pc,2pc>\xymatrix{
1\ar@{-}[dr]&&6\ar@{-}[dl]\ar@{-}[dr]&&8\ar@{-}[dl]\ar@{-}[dr]
&\save+<0.7pc,-1.7pc>\drop{\cdots}&&8\ar@{-}[dl]\ar@{-}[dr]\\ &3&&6&&6&6&&6
}\endxy$$
}
  
\noindent all of which are $\pinf (\lamod)$-phantoms of $S_1$ of the first
kind.  This yields the $\pinf (\lamod)$-phantom
$\varinjlim C_n$ of infinite $K$-dimension, and
shows that $S_1$ fails to  have a $\pinf
(\lamod)$-ap\-prox\-i\-ma\-tion inside $\pinf (\lamod).$ \qed\endexample

\example{Problem 14} Characterize the simple modules over monomial relation 
algebras which fail to have right $\pinf(\lamod)$-approximations 
in terms of their infinite dimensional
phantoms. \endexample

\head 4. Contravariant finiteness of $\pinf
(\lamod)$ and the inequality $\operatorname{fin\,dim} \Lambda <
\operatorname{Fin\,dim}
\Lambda$\endhead

We apply Criterion 10 to a less elementary example which, in fact, motivated a
major portion of this article. Namely, we show that, for the finite 
dimensional monomial relation algebra $\Lambda$ of [8] with $\lfindim \Lambda <
\lFindim \Lambda$, the category 
$\pinf (\lamod)$ is not contravariantly finite. 

\example{Example 15} We refer the reader to [8, p\. 378] for a
definition of $\Lambda = K\Gamma / I$. We will apply Criterion 10 to
show that the simple module $S_2 = \la e_2 / Je_2$ fails to have a right
$\pinf(\lamod)$-approximation.  For that purpose, we let $\frak A = \pinf(\lamod)$, set $m=1$, and
focus on the single primitive idempotent $e_2$.  Moreover, we make the choices $p = \gamma_1$
and
$q = \gamma_2 + \tau \gamma_2$, let $n \in \sN$, and set $x_i = e_2$ for
$i = 1,
\dots, n$. 

First we exhibit modules $M_n$ as in part (1) of the
criterion.  Namely, we define 
$$M_n = \left( \bigoplus^n_{i=1} \Lambda x_i \right) \bigg/ \left(
\sum^{n-1}_{i=1}
\Lambda z_i \right),$$
where $z_i = p x_i - q x_{i+1}$ for $1 \le i \le n-1$.  
Observe that $M_n \in \pinf(\lamod)$ for each $n$. Indeed, the sum
$\sum^{n-1}_{i=1}\Lambda z_i$ is direct and can be seen to have finite projective dimension as
follows:  The graph of
$\Lambda z_i$ relative to the top element $z_i$ is
  
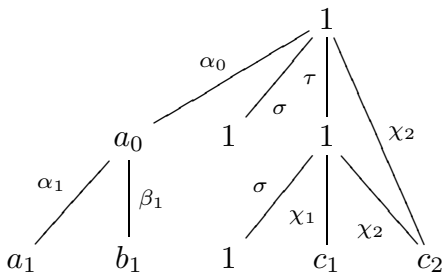
\ignore{
$$\xy\xymatrixcolsep{1.5pc}
\xymatrix{
&&&1 \ar@{-}[lld]_{\alpha_0}\ar@{-}[ld]^(0.65)\sigma\ar@{-}[d]_\tau
\ar@{-}[rdd]^{\chi_2}\\
&a_0\ar@{-}[dl]_{\alpha_1}\ar@{-}[d]^{\beta_1}
&1&1\ar@{-}[dl]_\sigma \ar@{-}[d]_(0.65){\chi_1}\ar@{-}[dr]_(0.65){\chi_2}\\ 
a_1&b_1&1&c_1&c_2
}\endxy$$  }
  
\noindent whence the graphical method of [7, Section 5] yields $\Omega^1 (\la
z_i) =
\ignore{
\xymatrixrowsep{.5pc}

\left( \xy<1pc,1pc>\xymatrix{c_1\ar@{-}[d]^{\chi_1'}\\
c_1}\endxy \right)\oplus
\left( \xy<1pc,1pc>\xymatrix{c_2\ar@{-}[d]^{\chi_2'}\\
c_2}\endxy \right)}$.  Thus 
$$\Omega^1 \left( \sum^{n-1}_{i=1}\Lambda z_i \right) \cong
\bigl(\la e(c_1)
\bigr)^{n-1} \oplus \bigl(\la e(c_2) \bigr)^{n-1}$$ 
is projective as required.

To check condition (2)(i) of Criterion 10, suppose that $C$ belongs to $\pinf
(\lamod)$ and has a top element $x$ of type $e_2$. Then $p x = \gamma_1 x \ne
0$, since otherwise
$\Omega^1 (C)$ would have a direct summand isomorphic to the left ideal $\la \gamma_1$.  But
this left ideal has infinite projective dimension, as can again be checked with the aid of its 
graph
  
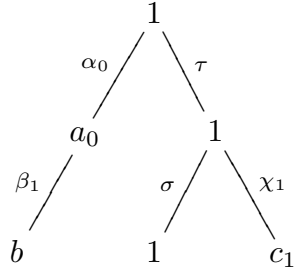
\ignore{
$$\xy\xymatrixcolsep{.7pc}
\xymatrix{
&&1\ar@{-}[dl]_{\alpha_0}\ar@{-}[dr]^\tau\\
&a_0\ar@{-}[dl]_{\beta_1}&&1\ar@{-}[dl]_\sigma\ar@{-}[dr]^{\chi_1}\\
b&&1&&c_1
}\endxy$$
}
  
\noindent and the method of [7].  This shows that condition (2)(i) is indeed
met.  

Finally, let us check that condition (2)(ii) of Criterion 10 is satisfied. 
Again let $C \in \pinf(\lamod)$, and suppose that
$y, z \in C$ are such that $0 \ne p y = q z$, where  $p$ and $q$ are as above.
From the fact that 
$qz = \gamma_2 z + \tau \gamma_2 z$ does not vanish,  we deduce $\gamma_2 z
\ne 0$, which in turn implies that $e_2 z$ is a top element of type $e_2$  of
$C$; this implication is an immediate consequence of the fact that the vertex
$e_2$ is a source of
$\Gamma$. Consequently, the preceding paragraph yields
$p z\ne 0$ as required.  Thus Criterion 10 applies to complete that proof that $S_2$ does
not have a right $\pinf(\lamod)$-approximation.

An infinite dimensional $\pinf (\lamod)$-phantom of $S_2$
resulting from the preceding argument can be visualized as follows:
  
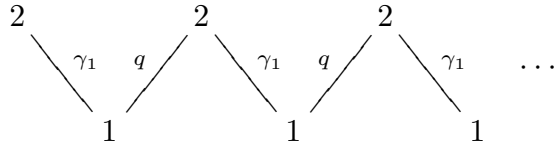
\ignore{
$$\xy\xymatrixcolsep{1.5pc}
\xymatrix{
2 \ar@{-}[dr]^{\gamma_1} &&2 \ar@{-}[dr]^{\gamma_1} \ar@{-}[dl]_q &&2
\ar@{-}[dr]^{\gamma_1} \ar@{-}[dl]_q \save+<4pc,-1.3pc>\drop{\cdots}&&\\
 &1 &&1 &&1
}\endxy$$
}
  
\noindent where again $q = \gamma_2 + \tau \gamma_2$. \qed \endexample

We believe that contravariant finiteness of $\pinf (\lamod)$ in
$\lamod$ is a condition strong enough to significantly impinge on the category
of arbitrary (not necessarily finitely generated)  
left $\la$-modules of finite projective
dimension.

\example{Problem 16} Decide whether contravariant finiteness of $\pinf (\lamod)$ 
implies equality
of the little and big left finitistic dimensions of $\la$. \endexample


\Refs
\widestnumber\key{HZ2}

\ref\no1\by J.L. Alperin \paper Diagrams for modules \jour J. Pure Appl.
Algebra \vol 16 \yr 1980 \pages 111-119 \endref 

\ref\no2\by M. Auslander, I. Reiten\paper Applications of
Contravariantly Finite Subcategories\jour Adv. in Math.\vol 86 \yr 1991\pages
111-152\endref

\ref\no3\by M. Auslander, S.O. Smal{\o}\paper Almost split sequences in
subcategories\jour J. Algebra\vol 69 \yr 1981\pages 426-454\endref

\ref\no4\by W.D. Burgess and B. Huisgen-Zimmermann \paper Approximating
modules by modules of finite projective dimension \jour J. Algebra \vol 178
\yr 1995 \pages 48-91
\endref

\ref\no5\by K.R. Fuller \paper Algebras from diagrams \jour J. Pure Appl.
Algebra \vol 48 \yr 1987 \pages 23-37 \endref
 
\ref\no6\by M. Harada, Y. Sai\paper On categories of indecomposable
modules I\jour Osaka J. Math.\vol 7 \yr 1970\pages 323-344\endref

\ref\no7\by B. Huisgen-Zimmermann\paper Predicting syzygies over
finite dimensional monomial relation algebras\jour Manuscr. Math.\vol 70 \yr
1991\pages 157-182\endref

\ref\no8\by B. Huisgen-Zimmermann\paper Homological domino effects and
the first Finitistic Dimension Conjecture\jour Invent. Math.\vol 108 \yr
1992\pages 369-383\endref

\ref\no9\by K. Igusa, S.O. Smal{\o}, G. Todorov\paper Finite
projectivity and contravariant finiteness\jour Proc. Amer. Math. Soc.\vol 109
\yr 1990\pages 937-941\endref

\endRefs
\enddocument